\documentstyle[12pt]{article}

\begin{document}
\author{S.V. L\"udkovsky}
\title{Non-Archimedean valued quasi-invariant descending
at infinity measures}
\date{10 May 2004}
\maketitle

\par This article is devoted to new results of investigations of
quasi-invariant non-Archimedean valued measures, which is becoming
more important nowdays due to the development of non-Archimedean
mathematical physics, particularly, quantum mechanics,
quantum field theory, theory of superstrings and supergravity
\cite{vla2,vla3,ardrvo,cas,djdr,khrum,lutmf99,yojan}.
On the other hand, quantum mechanics is based
on measure theory and probability theory.
For comparison references are given below also on works, where
realvalued measures on non-Archimedean spaces were studied.
Stochastic approach in quantum field theory is actively used
and investigated especially in recent years \cite{alkar,khrif,khipr,khrum}.
As it is well-known in the theory of functions great role is played
by continuous functions and differentiable functions.
\par In the classical
measure theory the analog of continuity is quasi-invariance relative to
shifts and actions of linear or non-linear operators in the Banach space,
differentiability of measures is the stronger condition and there is
very large theory about it in the classical case. Apart from it the
non-Archimedean  case was less studied. Since there are not
differentiable functions from the  field $\bf Q_p$ into $\bf R$
or in another non-Archimedean field ${\bf Q}_{p'}$ with $p\ne p'$,
then instead of differentiability of measures their pseudo-differentiability
is considered.
\par Effective ways to use quasi-invariant and pseudo-differentiable measures
are given in the articles of the author
\cite{luanma,luumn582,lujms112,lulapm,lubp99,lubp2,lutmf99,lufpmsp,
lufpm,luumn995}. I.V. Volovich was discussing with me the matter and
interested in results of my investigations
of non-Archimedean analogs of Gaussian measures such as to satisfy
as many Gaussian properties as possible as he has planned
to use such measures in non-Archimedean quantum field theory.
The question was not
so simple. He has supposed that properties with mean values,
moments, projections, distributions and convolutions of such measures
can be considered analogously. But thorough analysis has shown,
that not all properties can be satisfied, because in such case
the linear space would have a structure of the $\bf R$-linear space.
Nevertheless, many of the properties it is possible
to satisfy in the non-Archimedean case also. Gaussian measures
are convenient to work in the classical case, but in the non-Archimedean
case they do not play so great role.
\par Strictly speaking no any nontrivial Gaussian measure exists in
the non-Archimedean case, but measures having few properties
analogous to that of Gaussian can be outlined. Supplying them with
definite properties depends on a subsequent task for which problems
they may be useful. Certainly if each projection $\mu _Y$ of a measure
$\mu $ on a finite dimensional subspace $Y$ over a field $\bf K$
is equivalent to the Haar measure $\lambda _Y$ on $Y$, then this is well
property. But in the classical case, as it is well-known,
such property does not imply that the measure $\mu $ is Gaussian,
since each measure $\nu _Y(dx)=f(x)\lambda _Y(dx)$ with
$f\in L^1(Y,\lambda _Y,{\bf R})$ is absolutely continuous relative
to the Lebesgue measure $\lambda _Y$ on $Y$ and this does not imply
Gaussian properties of moments or its characteristic functional
\cite{gevi,dal}. The class of measures having such properties
of projections is described by the Kolmogorov and Kakutani theorems.
At first it is mentioned below how measures on Banach spaces can be used
for construction of measures on complete ultrauniform spaces,
then particular classes of quasi-invariant non-Archimedean valued measures
descending at infinity are considered.
\par In \cite{lufpm,luumn995} non-Archimedean polyhedral
expansions of ultrauniform spaces were investigated and the following
theorem was proved.
\par {\bf Theorem.} {\it Let $X$ be a complete ultrauniform space
and $\bf K$ be a local field. Then there exists an irreducible normal
expansion of $X$ into the limit of the inverse system
$S= \{ P_n,f^m_n,E \} $ of uniform polyhedra over $\bf K$, moreover,
$\lim S$  is uniformly isomorphic with $X$, where $E$ is an ordered
set, $f^m_n: P_m\to P_n$ is a continuous maping for each $m\ge n$;
particularly for the ultrametric space $(X,d)$ with the ultrametric
$d$ the inverse system $S$ is the inverse sequence.}
\par This structure theorem serves to prove the following theorem.
\par {\bf 1. Theorem.} {\it Let $X$ be a complete separable
ultrauniform space and let $\bf K$ be a local field.
Then for each marked $b\in \bf C_s$ there exists a nontrivial $\bf F$-valued
measure $\mu $ on $X$ which is a restriction of a measure $\nu $
in a measure space $(Y, Bco(Y),\nu ) = \lim \{ (Y_m, Bco(Y_m), \nu _m),
{\bar f}^m_n, E \} $ on $X$ and each $\nu _m$ is quasi-invariant and
pseudo-differentiable for $b\in \bf C_s$ relative to a dense subspace
${Y'}_m$, where $Y_n:=c_0({\bf K},\alpha _n)$,
${\bar f}^m_n: Y_m\to Y_n$ is a normal (that is,
$\bf K$-simplicial nonexpanding) mapping for each $m\ge n\in E$,
${\bar f}^m_n|_{P_m}=f^m_n$. Moreover, if $X$ is not locally
compact, then the family $\cal F$ of all such $\mu $ contains
a subfamily $\cal G$ of pairwise orthogonal measures
with the cardinality $card ({\cal G})=card ({\bf F})^{\sf c}$,
${\sf c}:=card ({\bf Q_p})$.}
\par {\bf Proof.} Choose a polyhedral expansion of $X$ in accordance
with cited above theorem. Let ${\bf Q_p}\subset \bf K$, $s\ne p$
are prime numbers, ${\bf Q_s}\subset \bf F$, where $\bf F$ is a
non-Archimedean field complete relative to its uniformity.
On each $X_n$ take a probability $\bf F$-valued measure $\nu _n$
such that $\| X_n\setminus P_n \| _{\nu _n}<\epsilon _n$, $\sum_{n\in E}
\epsilon _n< 1/5 $. In accordance with \S 3.5.1 and \S 4.2.1
\cite{lulapm,lujms112} (see also \cite{luumn582}) each $\nu _n$
can be chosen quasi-invariant and pseudo-differentiable for $b\in \bf C_s$
relative to a dense $\bf K$-linear subspace ${Y'}_n$, since each
normal mapping $f^m_n$ has a normal extension on $Y_m$ supplied
with the uniform polyhedra structure.
Since $E$ is countable and ordered, then a family $\nu _n$ can be
chosen by transfinite induction consistent, that is,
${\bar f}^m_n(\nu _m)=\nu _n$ for each $m\ge n$ in $E$,
${\bar f}^m_n({Y'}_m)={Y'}_n$. Then $X=\lim \{ P_m,f^m_n,E \}
\hookrightarrow Y$. Since ${\bar f}^m_n$ are $\bf K$-linear, then
$({\bar f}^m_n)^{-1}(Bco (Y_n))\subset Bco (Y_m)$ for each
$m\ge n\in E$. Therefore, $\nu $ is correctly defined
on the algebra $\bigcup_{n\in E}f_n^{-1} (Bco (Y_n))$ of subsets of $Y$,
where $f_n: X \to X_n$ are $\bf K$-linear continuous epimorphisms.
Since $\nu $ is nontrivial and $ \| \nu \| $ is bounded by $1$, then by the
non-Archimedean analog of the Kolmogorov theorem \cite{lufpmsp,lukhr}
$\nu $ has an
extension on the algebra $Bco (Y)$ and hence on its completion $Af (Y,\nu )$.
Put $Y':=\lim \{ {Y'}_m, {\bar f}^m_n, E \} $. Then $\nu _m$ on $Y_m$
is quasi-invariant and pseudo-differentiable for $b\in \bf C_s$
relative to ${Y'}_m$. From $\sum_n\epsilon _n<1/5$ it follows, that
$1\ge \| X \| _{\mu }\ge \prod_n (1-\epsilon _n) > 1/2$, hence $\mu $ is
nontrivial. 
\par To prove the latter statement use the non-Archimedean analog
of the Kakutani theorem (see \cite{lulapm,lujms112})
for $\prod_nY_n$ and then consider the embeddings
$X\hookrightarrow Y\hookrightarrow \prod_nY_n$ such that projection
and subsequent restriction of the measure $\prod_n\nu _n$ on $Y$
and $X$ are nontrivial, which is possible due to the proof given above.
If $\prod_n\nu _n$ and $\prod_n{\nu '}_n$ are orthogonal on $\prod_nY_n$,
then they give $\nu $ and $\nu '$ orthogonal on $X$.
\par {\bf 2. Definitions and Notes.} A function $f: {\bf K}\to 
\bf U_s$ is called pseudo-differentiable of order $b$, if there exists
the following integral:
$PD(b,f(x)):=\int_{\bf K}[(f(x)-f(y)) \times g(x,y,b)]
dv(y)$. We introduce the following notation $PD_c(b,f(x))$ for
such integral by $B({\bf K},0,1)$ instead of the entire $\bf K$.
Where $g(x,y,b):=s^{(-1-b) \times ord_p(x-y)}$
with the corresponding Haar measure $v$
with values in $\bf K_s$, where $\bf K_s$ is a local field containing
the field $\bf Q_s$, $s$ is a prime number, $b \in {\bf C_s}$ and
$|x|_{\bf K}=p^{-ord_p(x)}$,
$\bf C_s$ denotes the field of complex numbers with the
non-Archimedean valuation extending that of $\bf Q_s$, 
$\bf U_s$ is a spherically complete field with a valuation group
$\Gamma _{\bf U_s} := \{ |x|:$ $0\ne x \in {\bf U_s} \} =
(0,\infty )\subset \bf R$ such that ${\bf C_s}\subset \bf U_s$,
$0<s$ is a prime number \cite{diar,roo,sch1,wei}.
For each $\gamma \in (0,\infty )$ there exists $\alpha =log_s(\gamma )
\in  \bf R$, $\Gamma _{\bf U_s}=(0,\infty )$, hence
$s^{\alpha }\in \bf U_s$ is defined for each $\alpha \in \bf R$,
where $log_s(\gamma )= ln (\gamma )/ ln (s) $,
$ln: (0,\infty )\to \bf R$ is the natural logarithmic
function such that $ln (e)=1$.
The function $s^{\alpha +i \beta } =:
\xi (\alpha ,\beta )$ with $\alpha $ and $\beta \in \bf R$
is defined due to the algebraic isomorphism of $\bf C_s$ with
$\bf C$ (see \cite{kobl}) in the following manner.
Put $s^{\alpha + i\beta }:=s^{\alpha }(s^i)^{\beta }$
and choose as $s^i$ a marked number in $\bf U_s$ such that
$s^i:=(EXP_s(i))^{ln~ s}$, where $EXP_s: {\bf C_s}\to \bf C_s^+$
is the exponential function, ${\bf C_s^+} := \{ x\in {\bf C_s}:
|x-1|_s<1 \} $ (see Proposition $45.6$ \cite{sch1}).
Therefore, $|EXP_s(i) -1|_s<1$, hence
$|EXP_s(i)|_s=1$ and inevitably $|s^i|_s=1$. Therefore,
$|s^{\alpha +i\beta }|_s=s^{-\alpha }$ for each $\alpha $ and $\beta
\in \bf R$, where $|*|_s$ is the extension of the valuation
from $\bf Q_s$ on $\bf U_s$, consequently, $s^x\in \bf U_s$
is defined for each $x\in \bf C_s$.
\par A quasi-invariant measure $\mu $ on $X$
is called pseudo-differentiable for $b \in \bf C_s$,
if there exists $PD(b,g(x))$ for $g(x):=\mu (-xz+S)$ for each
$S \in Bco(X)$ $\| S \| _{\mu }< \infty $
and each $z \in J^b_{\mu }$, where $J^b_{\mu }$ is a $\bf K$-linear
subspace dense in $X$. For a fixed $z \in X$ such measure is called
pseudo-differentiable along $z$.
\par {\bf 2.1. Definitions and Remarks.} Let $X$ be a locally
$\bf K$-convex space equal to
a projective limit $\lim \{ X_j, \phi ^j_l, \Upsilon \} $
of Banach spaces over a local field $\bf K$ such that
$X_j=c_0(\alpha _j,{\bf K})$, where the latter space consists
of vectors $x=(x_k: k\in \alpha _j )$, $x_k\in \bf K$,
$\| x \| :=\sup_k |x_k|_{\bf K}<\infty $ and such that
for each $\epsilon >0$ the set $\{ k: |x_k|_{\bf K}>\epsilon \} $
is finite, $\alpha _j$ is a set, that is convenient
to consider as an ordinal due to Kuratowski-Zorn lemma \cite{eng,roo};
$\Upsilon $ is an ordered set,
$\phi ^j_l: X_j\to X_l$ is a $\bf K$-linear continuous mapping
for each $j\ge l\in \Upsilon $, $\phi _j: X\to X_j$ is a projection
on $X_j$, $\phi _l\circ \phi ^j_l=\phi _j$ for each $j\ge l\in \Upsilon $,
$\phi ^l_k\circ \phi ^j_l=\phi ^j_k$ for each $j\ge l\ge k$ in $\Upsilon $. 
Consider also a locally $\bf R$-convex space, that is
a projective limit $Y=\lim \{ l_2(\alpha _j,{\bf R}),
\psi ^j_l, \Upsilon \} $, where $l_2(\alpha _j,{\bf R})$ is the real Hilbert
space of the topological weight $w(l_2(\alpha _j,{\bf R}))=
card (\alpha _j)\aleph _0$. Suppose $B$ is a symmetric
nonegative definite (bilinear) nonzero functional $B: Y^2\to \bf R$.
\par Consider a non-Archimedean field $\bf F$
such that ${\bf K_s} \subset \bf F$ and with the valuation group
$\Gamma _{\bf F}=(0,\infty )\subset \bf R$ and $\bf F$ is complete
relative to its uniformity (see \cite{diar,esc}).
Then a measure $\mu =\mu _{q,B,\gamma }$ on $X$
with values in $\bf K_s$ is called a $q$-Gaussian measure, if
its characteristic functional $\hat \mu $ with values in $\bf F$ has the form
$${\hat \mu } (z)=s^{[B(v^s_q(z),v^s_q(z))]}\chi _{\gamma }(z)$$
on a dense $\bf K$-linear subspace ${\sf D}_{q,B,X}$ in $X^*$ of all
continuous $\bf  K$-linear functionals $z: X\to \bf K$
of the form $z(x)=z_j(\phi _j(x))$ for each $x\in X$
with $v^s_q(z)\in {\sf D}_{B,Y}$,
where $B$ is a nonnegative definite bilinear $\bf R$-valued
symmetric functional on a dense $\bf R$-linear subspace
${\sf D}_{B,Y}$ in $Y^*$, $B: {\sf D}_{B,Y}^2\to \bf R$,
$j\in \Upsilon $ may depend on $z$, $z_j: X_j\to \bf K$ is a continuous
$\bf K$-linear functional such that $z_j=\sum_{k\in \alpha _j}
e^k_jz_{k,j}$ is a countable convergent series such that
$z_{k,j}\in \bf K$, $e^k_j$ is a continuous $\bf K$-linear functional
on $X_j$ such that $e^k_j(e_{l,j})=\delta ^k_l$ is the Kroneker delta
symbol, $e_{l,j}$ is the standard orthonormal (in the non-Archimedean sence)
basis in $c_0(\alpha _j,{\bf K})$, $v^s_q(z)=v^s_q(z_j):=\{
|s^{q~ ord_p(z_{k,j})/2}|_s: k\in \alpha _j \} $.
It is supposed that $z$ is such that $v^s_q(z)\in l_2(\alpha _j,{\bf R})$,
where $q$ is a positive constant, $\chi _{\gamma }(z): X\to \bf T_s$
is a continuous character such that $\chi _{\gamma }(z)=
\chi (z(\gamma ))$, $\gamma \in X$, $\chi : {\bf K}\to \bf T_s$
is a nontrivial character of $\bf K$ as an additive group
(see \cite{roo} and \S 2.5 in \cite{lulapm,lujms112}).
\par {\bf 3. Proposition.} {\it A $q$-Gaussian quasi-measure
on an algebra of cylindrical subsets $\bigcup_j\pi _j^{-1}({\cal R}_j)$,
where $X_j$ are finite-dimensional over $\bf K$ subspaces in $X$, is a
measure on a covering ring ${\cal R}$ of subsets of $X$ (see \S 2.36
\cite{lulapm,lujms112}). Moreover, a correlation operator $B$ is of
class $L_1$, that is, $Tr (B)<\infty $, if and only if
each finite dimensional over $\bf K$
projection of $\mu $ is a $q$-Gaussian measure (see \S 2.1).}
\par {\bf Proof.} From Definition 2.1 it follows, that
each one dimensional over $\bf K$ projection $\mu _{x\bf K}$
of a measure $\mu $ satisfies Conditions $2.1.(i-iii)$
\cite{lulapm,lujms112}
the covering ring $Bco ({\bf K})$, where $0\ne x=e_{k,l}\in X_l$.
Therefore, $\mu $ is defined and finite additive on a cylindrical algebra
\par ${\sf U}:=\bigcup_{k_1,...,k_n;l}\phi _l^{-1}
[(\phi ^l_{k_1,...,k_n})^{-1}
(Bco (span_{\bf K}\{ e_{k_1,l},...,e_{k_n,l} \} ))]$, \\
where $\phi ^l_{k_1,...,k_n}: X_l\to span_{\bf K}(e_{k_1,l},...,e_{k_n,l})$
is a projection. This means that $\mu $ is a bounded quasimeasure
on $\sf U$. Since ${\hat {\mu }}(0)=1$, then $\mu (X)=1$.
The characteristic functional ${\hat {\mu }}$ satisfies Conditions
$2.5.(3,5)$ \cite{lulapm,lujms112}. In view of the non-Archimedean
analog of the Bochner-Kolmogorov theorem \S 2.21
and Theorem $2.37$ \cite{lulapm,lujms112}
$\mu $ has an extension to a probability measure on
a covering ring $\cal R$ of subsets of $X$ containing $\sf U$.
\par Suppose that $B$ is of class $L_1$.
Then $B(v_q(z),v_q(z))$ and hence $\hat \mu (z)$ is correctly
defined for each $z\in {\sf D}_{q,B,X}$.
The set ${\sf D}_{q,B,X}$ of functionals $z$ on $X$ from
\S 2.1 separates points of $X$.
From Definition 2.1 it follows, that ${\hat {\mu }}(y)$
is continuous. Consider a diagonal compact operator 
$T$ in the standard orthonormal base, $Te_{k,l}=a_{k,l}e_{k,l}$,
$\lim_{k+l\to \infty }a_{k,l}=0$. Since $B$ is continuous, then
the corresponding to $B$ correlation operator $E$ is a bounded
$\bf K$-linear operator on $Y$, $\| E \| <\infty $.
For each $\epsilon >0$ there exist $\delta >0$ 
and $T$ such that $\max (1, \| E \| ) \delta <\epsilon $ and
$|a_{k,l}|<\delta $ for each $k+l>N$, where $N$ is a marked natural number,
therefore, $\| E|_{span_{\bf K} \{ e_{k,l}: k+l>N \} } \| <
\epsilon $. Hence for each $\epsilon >0$ there exists a compact operator
$T$ such that from $|{\tilde z}Tz|<1 $ it follows,
$|{\hat {\mu }}(y)-{\hat {\mu }}(x)|< \epsilon $ for each $x-y=z$,
where $x, y, z\in Y^*$. Therefore, by Theorem $2.30$ the charateristic
functional ${\hat {\mu }}$ defines a probability Radon measure
on $Bco (X)$.
\par Vice versa suppose that each finite dimensional over $\bf K$
projection of $\mu $ is a measure of the same type.
If for a given one dimensional over $\bf K$ subspace $W$
in $X$ it is the equality $B(v_q(z),v_q(z))=0$ for each
$z\in W$, then the projection $\mu _W$ of $\mu $
is the atomic measure with one atom. Show
$B\in L_1(c_0(\omega _0,{\bf K}))$ and $\gamma \in c_0(\omega _0,{\bf K})$.
Let $0\ne x\in X$ and consider the projection $\pi _x: X\to x\bf K$.
Since $\mu _{x\bf K}$ is the measure on $Bco (x{\bf K})$,
then its characterisic functional satisfies Conditions of Theorem $2.30$
\cite{lulapm,lujms112}.
Then ${\hat {\mu }}$ for $x\bf K$ gives the same characteristic functional
of the type
$${\hat {\mu }}_{x\bf K}(z)=s^{[b_x(v^s_q(z))^2]}\chi _{{\delta }_x}(z)$$
for each $z\in x\bf K$, where $b_x>0$ and $\delta _x\in \bf K$ are constants
depending on the parameter $0\ne x\in X$. Since $x$ and $z$ are arbitrary,
then this implies, that $B\in L_1$ and $\gamma \in c_0(\omega _0,{\bf K})$.
\par {\bf 4. Corollary.} {\it A $q$-Gaussian measure $\mu $
from Proposition 3 with $Tr (B)<\infty $ is quasi-invariant
and pseudo-differentiable for some $b\in \bf C_s$
relative to a dense subspace $J_{\mu }\subset M_{\mu }=
\{ x\in X:$ $v_q^s(x)\in E^{1/2}(Y) \} $.
Moreover, if $B$ is diagonal, then each one-dimensional projection $\mu ^g$
has the following characteristic functional:
$$(i)\quad {\hat \mu }^g(h)=s^{(\sum_j\beta _j|g_j|^q)|h|^q}
\chi _{g(\gamma )}(h),$$
where $g=(g_j: j\in \omega _0)\in c_0(\omega _0,{\bf K})^*$,
$\beta _j>0$ for each $j$.}
\par {\bf Proof.} Using the projective limit reduce consideration to
the Banach space $X$. Take a prime number $s$ such that $s\ne p$
and consider a field $\bf K_s$ such that $\bf K$ is compatible
with $\bf K_s$, which is possible, since $\bf K$ is a finite
algebraic extension of $\bf Q_p$ and it is possible to take
in particular ${\bf K_s}=\bf Q_s$. Recall that a group $G$ for which
$o(G)\subset o({\bf T}_{\bf K})$ is called compatible with $\bf K$,
where $o(G)$ denotes the set of all natural numbers for which $G$ has
an open subgroup $U$ such that at least one of the elements
of the quotient group $G/U$ has order $n$, $\bf T$ denotes the group of all
roots of $1$ and ${\bf T}_{\bf K}$ denotes its subgroup of all elements
whose orders are not divisible by the characteristic $p$ of the residue
class field $k$ of $\bf K$. A character of $G$ is a continuous
homomorphism $f: G\to \bf T$. Under pointwise multiplication charaters
form a group denoted by $G^{\hat .}$.
A group $G$ is called torsional, if each compact subset $V$ of $G$
is contained in a compact subgroup of $G$.
In view of Theorem $9.14$ \cite{roo}
${\bf K}^{\hat .}$ is isomorphic with $\bf K$. A $\bf K$-valued
character of a group $G$ is a continuous homomorphism
$f: G\to {\bf T}_{\bf K}$.
The family of all $\bf K$-valued characters form a group denoted by
$G^{\hat .}_{\bf K}$. Since $\bf K$ is compatible with $\bf K_s$
and $\lim_{n\to \infty }p^n=0$, then ${\bf K}^{\hat .}$ is isomorphic
with ${\bf K}^{\hat .}_{\bf K_s}$.
If $G$ is a torsional group, then the Fourier-Stieltjes transform
of a tight measure $\mu \in M(G)$ is the mapping ${\hat {\mu }}:
G^{\hat .}_{\bf K}\to \bf K$ defined by the formula:
${\hat {\mu }}(g):=\int_G \chi (x)\mu (dx)$, where $\chi \in
G^{\hat .}_{\bf K}$. In view of Schikhof Theorem $9.21$
\cite{roo} the Fourier-Stieltjes transformation induces a Banach
algebra isomorphism $L(G,{\cal R},w,{\bf K})$ with
$C_{\infty }(G^{\hat .}_{\bf K},{\bf K})$, where $w$ is
a nontrivial Haar $\bf K$-valued measure on $G$. Therefore,
in this sutuation there exists the  Banach algebra isomorphism
of $L({\bf K},{\cal R},w,{\bf K_s})$ with
$C_{\infty }({\bf K}^{\hat .}_{\bf K_s},{\bf K_s})$.
\par Therefore, from the proof above and Theorem $3.5$ it follows, 
that the measure $\mu _{q,B,\gamma }$ is quasi-invariant relative to shifts
on vectors from the dense subspace $X'$ in $X$ such that
$X'=\{ x\in X: v^s_q(x) \in E^{1/2}(Y) \} $, which is $\bf K$-linear,
since $B$ is $\bf R$-bilinear and $B(y,z)=:(Ey,z)$ for each $y, z\in Y$
and $v^s_q(ax)=|a|^{q/2}v^s_q(x)$ and $v^s_q(x_j+t_j)\le
\max (v^s_q(x_j),v^s_q(t_j))$ for each $x, t\in X$ and each $a\in \bf K$,
where $E$ is nondegenerate positive definite of trace class
$\bf R$-linear operator on $Y$, $x=\sum_jx_je_j$, $x_j\in \bf K$,
since $l_2^*=l_2$ and $E$ can be
extended from ${\sf D}_{B,Y}$ on $Y$.
\par Consider $s^{a+ib}$ as in \S 2. Mention, that $|(|z|_p)|_s=1$ for
each $z\in \bf K$, where the field $\bf K$ is compatible with $\bf K_s$.
\par The pseudo-differential operator has the form:
$PD(b,f(x)):=\int_{\bf K}[f(x)-f(y)] s^{(-1-b)\times ord_p(x-y)} w(dy)$,
where $w$ is the Haar $\bf K_s$-valued measure on $Bco ({\bf K})$,
$b\in \bf C_s$, particularly, also for $f(x):=\mu (-xz+A)$
for a given $z\in X'$, $A\in Bco (X)$, where $x, y\in \bf K$.
Using the Fourier-Stieltjes transform write it in the form:
$PD(b,f(x))=F^{-1}_v(\xi (v) \psi (v))$,
where $\xi (v) := [F_y(f(x)-f(y))](v)$, $\psi (v):=
[F_y (s^{(-1-b)\times ord_p(y)})] (v)$, $F_y$ means the Fourier-Stieltjes
operator by the variable $y$. Denoting $A-xz=:S$ we can consider
$f(x)=0$ and $f(y)=\mu ((x-y) z + S) - \mu (S)$, since $S\in Bco (X)$.
Then $f(y)=\int_S(\mu ((x-y)+dg)-\mu (dg))=\int_S[\rho _{\mu }(y-x,g)
-1]\mu (dg)$. The constant function $h(g)=1$ is evidently
pseudo-differentiable of order $b$ for each $b\in \bf C_s$.
Hence the pseudo-differentiability of $\mu $ of order $b$
follows from the existence of pseudo-differential of
the quasi-invariance factor $\rho _{\mu }(y,g+x)$ of order $b$
for $\mu $-almost every $g\in X$. In view of Theorem $3.5$
and the Fourier-Stieltjes operator isomorphism of Banach algebras
$L({\bf K},{\cal R},w,{\bf K_s})$ and
$C_{\infty }({\bf K}^{\hat .}_{\bf K_s},{\bf K_s})$ the
pseudo-differentiability of $\rho _{\mu }$ follows from the existence
of $F^{-1}({\hat {\mu }} \psi )$, where ${\hat {\mu }}$ is the characteristic
functional of $\mu $. We have \\
$(ii)\quad F(f)(y)=\int_{\bf K} \chi (xy)f(x)w(dx)$ \\
$=\int_{\bf K}\chi (z) f(z/y)[|y|_p]^{-1}w(dz)$ \\
for each $y\ne 0$, where $x, y, z\in \bf K$, particularly, for $f(x)=
s^{-(1+b)\times ord_p(x)}$ we have $f(z/y)=f(z) f(-y)$ and
$F(f)(y) = \Gamma ^{{\bf K},s}(1+b)f(-y)|y|_p^{-1}$, where \\
$(iii)\quad \Gamma ^{{\bf K},s}(b):=
\int_{\bf K}\chi (z)s^{-b \times ord_p(x)}w(dz)$, \\
$f(-y)=s^{(1+b)\times ord_p(y)}$, since $ord_p(z/y)=ord_p(x) - ord_p(y)$.
For a nontrivial character of an order $m\in \bf Z$ from the definition
it follows, that $\Gamma ^{{\bf K},s}(b)\ne 0$
for each $b$ with $Re (b)\ne 0$,
since $|s^{-bn}|_s=s^{Re (b)n}$ for each $n\in \bf Z$.
Therefore, $\psi (y)=s^{(1+b)\times ord_p(y)}|y|_p^{-1}$, consequently,
$|\psi (y)|_s = s^{-(1+Re (b))\times ord_p(y))}$ for each $y\ne 0$,
since $|(|y|_p)|_s=1$. On the other hand, $|{\hat {\mu }}(z)|=
s^{-B(v^s_q(z),v^s_q(z))}$ and $F^{-1}({\hat {\mu }} \psi )$
exists for each $b\in \bf C_s$ with $Re (b)>-1$, since $Tr (B)<\infty $,
which is correct, since
$\bf C_s$ is algebraically isomorphic with $\bf C$ and $\Gamma _{\bf U_s}
\supset (0,\infty )$.
\par {\bf 5. Corollary.} {\it Let $X$ be a complete locally
$\bf K$-convex space of separable type over a local field $\bf K$,
then for each constant $q>0$ there exists a nondegenerate
symmetric positive definite operator $B\in L_1$ such that
a $q$-Gaussian quasi-measure is a measure on $Bco (X)$ and
each its one dimensional over $\bf K$ projection is absolutely
continuous relative to the nonnegative Haar measure on $\bf K$.}
\par {\bf Proof.} A space $Y$ from \S 2.1 corresponding to
$X$ is a separable locally $\bf R$-convex space.
Therefore, $Y$ in a weak topology is isomorphic with
${\bf R}^{\aleph _0}$ from which the existence of $B$ follows.
For each $\bf K$-linear finite dimensional over $\bf K$ subspace
$S$ a projection $\mu ^S$ of $\mu $ on $S\subset X$
exists and its density $\mu ^S(dx)/w(dx)$ relative to the
nondegenerate $\bf K_s$-valued Haar measure $w$ on $S$ is the inverse
Fourier-Stieltjes transform $F^{-1}({\hat \mu }|_{S^*})$
of the restriction of ${\hat \mu }$ on $S^*$.
For $B\in L_1$ each one dimensional projection of $\mu $
corresponding to $\hat \mu $ has a density that is a
continuous function belonging to $L({\bf K},Bco ({\bf K}),w,{\bf K_s})$.
\par {\bf 6. Proposition.} {\it Let $\mu _{q,B,\gamma }$ and
$\mu _{q,E,\delta }$ be two $q$-Gaussian measures with
correlation operators $B$ and $E$ of class $L_1$,
then there exists a convolution
of these measures $\mu _{q,B,\gamma }*\mu _{q,E,\delta }$, which is
a $q$-Gaussian measure $\mu _{q,B+E,\gamma +\delta }$.}
\par {\bf Proof.} Since $B$ and $E$ are nonnegative, then
$(B+E)(y,y)=B(y,y)+E(y,y)\ge 0$ for each $y\in Y$, that is,
$B+E$ is nonnegative. Evidently, $B+E$ is symmetric
and of class $L_1$. Moreover, $\mu _{q,B+E,\gamma +\delta }$ is
defined on the covering ring ${\sf U}_{B+E}$ containing the union of
covering rings ${\sf U}_B$ and ${\sf U}_E$
on which $\mu _{q,B,\gamma }$ and $\mu _{q,E,\delta }$
are defined correspondingly, since $ker (B+E)\subset ker (B)\cap ker (E)$.
Therefore, $\mu _{q,B+E,\gamma +\delta }$ is the tight
$q$-Gaussian measure together with $\mu _{q,B,\gamma }$ and
$\mu _{q,E,\delta }$ in accordance with Proposition 3
on the covering ring ${\cal R}_{\mu _{q,B+E,\gamma +\delta }}$
which is the completion of the minimal ring generated by ${\sf U}_{B+E}$.
Since ${\hat \mu }_{q,B+E,\gamma +\delta }={\hat \mu }_{q,B,\gamma }
{\hat \mu }_{q,E,\delta }$,
then $\mu _{q,B+E,\gamma +\delta }=\mu _{q,B,\gamma } *\mu _{q,E,\delta }$.
\par {\bf 6.1. Remark and Definition.} A measurable space $(\Omega ,{\sf F})$ 
with a probability $\bf K_s$-valued
measure $\lambda $ on a covering ring $\sf F$ of a set
$\Omega $ is called a probability space and it is denoted by
$(\Omega ,{\sf F},\lambda )$. Points $\omega \in \Omega $ are called 
elementary events and values $\lambda (S)$  
probabilities of events $S\in \sf F$. A measurable map
$\xi : (\Omega ,{\sf F})\to (X,{\sf B})$ is called a random variable
with values in $X$, where ${\sf B}$ is a covering ring such that
${\sf B}\subset Bco(X)$, $Bco(X)$ is the ring
of all clopen subsets of a locally $\bf K$-convex space $X$,
$\xi ^{-1}({\sf B})\subset \sf F$,
where $\bf K$ is a non-Archimedean field complete as an ultrametric space. 
\par The random variable $\xi $ induces a normalized measure
$\nu _{\xi }(A):=
\lambda (\xi ^{-1}(A))$ in $X$ and a new probability space 
$(X,{\sf B},\nu _{\xi }).$
\par Let $T$ be a set with a covering ring $\cal R$ and a measure
$\eta : {\cal R}\to \bf K_s$. Consider the following Banach space
$L^q(T,{\cal R},\eta ,H)$ as the completion of the set
of all ${\cal R}$-step functions $f: T\to H$ relative to the following 
norm:
\par $(1)\quad \| f\|_{\eta ,q}:=\sup_{t\in T}\| f(t)\|_H
N_{\eta }(t)^{1/q}$ for $1\le q<\infty $ and
\par $(2)\quad \| f\|_{\eta ,\infty }:=\sup_{1\le q<\infty }
\| f(t)\|_{\eta ,q}$, where $H$ is a Banach space over $\bf K$.
For $0<q<1$ this is the metric space with the metric
\par $(3)\quad \rho _q(f,g):=\sup_{t\in T}\| f(t)-g(t)\|_H
N_{\eta }(t)^{1/q}.$  
\par If $H$ is a complete locally $\bf K$-convex space,
then $H$ is a projective limit of Banach spaces
$H=\lim \{ H_{\alpha },\pi ^{\alpha }_{\beta }, \Upsilon \} $,
where $\Upsilon  $ is a directed set, $\pi ^{\alpha }_{\beta }:
H_{\alpha }\to H_{\beta }$ is a $\bf K$-linear continuous mapping
for each $\alpha \ge \beta $, $\pi _{\alpha }: H\to H_{\alpha }$
is a $\bf K$-linear continuous mapping such that $\pi ^{\alpha }_{\beta }
\circ \pi _{\alpha }=\pi _{\beta }$ for each $\alpha \ge \beta $
(see \S 6.205 \cite{nari}).
Each norm $p_{\alpha }$ on $H_{\alpha }$ induces a prednorm
${\tilde p}_{\alpha }$ on $H$. If $f: T\to H$, then $\pi _{\alpha }\circ
f=:f_{\alpha }: T\to H_{\alpha }$. In this case $L^q(T,{\cal R},\eta ,H)$
is defined as a completion of a family of all step functions
$f: T\to H$ relative to the family of prednorms
\par $(1')\quad \| f\|_{\eta ,q,\alpha }:=\sup_{t\in T}{\tilde p}_{\alpha }
(f(t))N_{\eta }(t)^{1/q}$, $\alpha \in \Upsilon $, for $1\le q<\infty $ and
\par $(2')\quad \| f\|_{\eta ,\infty ,\alpha  }:=\sup_{1\le q<\infty }
\| f(t)\|_{\eta ,q,\alpha }$, $\alpha \in \Upsilon $,
or pseudometrics
\par $(3')\quad \rho _{q,\alpha }(f,g):=\sup_{t\in T}{\tilde p}_{\alpha }
(f(t)-g(t))N_{\eta }(t)^{1/q}$, $\alpha \in \Upsilon $, for $0<q<1$.
Therefore, $L^q(T,{\cal R},\eta ,H)$ is isomorphic with the projective limit \\
$\lim \{ L^q(T,{\cal R},\eta ,H_{\alpha }),\pi ^{\alpha }_{\beta },
\Upsilon \} $.
For $q=1$ we write simply $L(T,{\cal R},\eta ,H)$ and
$\| f\|_{\eta }$. This definition is correct, since
$\lim_{q\to \infty }a^{1/q}=1$ for each $\infty >a>0$.
For example, $T$ may be a subset of $\bf R$. Let $\bf R_d$
be the field $\bf R$ supplied with the discrete topology. Since
the cardinality $card ({\bf R})={\sf c}=2^{\aleph _0},$
then there are bijective mappings of $\bf R$ on $Y_1:=\{ 0,...,b \}^{\bf N}$
and also on $Y_2:={\bf N}^{\bf N}$, where $b$ is a positive integer number.
Supply $\{ 0,...,b \}$ and $\bf N$ with the discrete topologies
and $Y_1$ and $Y_2$ with the product topologies.
Then zero-dimensional spaces
$Y_1$ and $Y_2$ supply $\bf R$ with covering separating rings
${\cal R}_1$ and ${\cal R}_2$ contained in $Bco(Y_1)$ and $Bco(Y_2)$
respectively. Certainly this is not related
with the standard (Euclidean) metric in $\bf R$.
Therefore, for the space $L^q(T,{\cal R},\eta ,H)$ we can consider
$t\in T$ as the real time parameter. If $T\subset \bf F$
with a non-Archimedean field $\bf F$, then we can consider
the non-Archimedean time parameter.
\par If $T$ is a zero-dimensional $T_1$-space, then denote by
$C^0_b(T,H)$ the Banach space of all continuous bounded functions
$f: T\to H$ supplied with the norm:
\par $(4)\quad \| f\|_{C^0}:=\sup_{t\in T} \| f(t)\|_H<\infty $. \\
If $T$ is compact, then $C^0_b(T,H)$ is isomorphic with
the space $C^0(T,H)$ of all continuous functions $f: T\to H$. 
\par For a set $T$ and a complete locally $\bf K$-convex 
space $H$ over $\bf K$ consider the 
product $\bf K$-convex space $H^T:=\prod_{t\in T}H_t$ in the product topology,
where $H_t:=H$ for each $t\in T$.
\par Then take on either $X:=X(T,H)=L^q(T,{\cal R},\eta ,H)$ or $X:=
X(T,H)=C^0_b(T,H)$ or on $X=X(T,H)=H^T$ a covering ring ${\sf B}$ such that
${\sf B}\subset Bco(X)$. Consider a random variable
$\xi : \omega \mapsto \xi (t,\omega )$ with values in $(X,{\sf B})$,
where $t\in T$.
\par Events $S_1,...,S_n$ are called independent in total if
$P(\prod_{k=1}^nS_k)=\prod_{k=1}^nP(S_k)$. Subrings
${\sf F}_k\subset {\sf F}$ are said to be independent if
all collections of events $S_k\in {\sf F}_k$ are independent in total, 
where $k=1,...,n$, $n\in \bf N$. To each collection of random variables
$\xi _{\gamma }$ on $(\Omega ,{\sf F})$ with $\gamma \in \Upsilon $
is related the minimal ring ${\sf F}_{\Upsilon }\subset \sf F$
with respect to which all $\xi _{\gamma }$ are measurable, where $\Upsilon $
is a set.
Collections $\{ \xi _{\gamma }: $ $\gamma \in \Upsilon _j \} $
are called independent if such are 
${\sf F}_{\Upsilon _j}$, where $\Upsilon _j\subset \Upsilon $ for each
$j=1,...,n,$ $n\in \bf N$.
\par Consider $T$ such that $card(T)>n$. For $X=C^0_b(T,H)$ or $X=H^T$
define $X(T,H;(t_1,...,t_n);(z_1,...,z_n))$ as a closed submanifold
in $X$ of all $f: T\to H$, $f\in X$ such that $f(t_1)=z_1,...,f(t_n)=z_n$, where
$t_1,...,t_n$ are pairwise distinct points in $T$ and
$z_1,...,z_n$ are points in $H$.
For $X=L^q(T,{\cal R},\eta ,H)$ and 
pairwise distinct points $t_1,...,t_n$ in $T$ with
$N_{\eta }(t_1)>0,...,N_{\eta }(t_n)>0$ define 
$X(T,H;(t_1,...,t_n);(z_1,...,z_n))$ as a closed submanifold
which is the completion relative to the norm $\| f\|_{\eta ,q}$
of a family of $\cal R$-step functions $f: T\to H$ such that
$f(t_1)=z_1,...,f(t_n)=z_n$. In these cases 
$X(T,H;(t_1,...,t_n);(0,...,0))$ is the proper $\bf K$-linear subspace 
of $X(T,H)$ such that $X(T,H)$ is isomorphic with
$X(T,H;(t_1,...,t_n);(0,...,0))\oplus H^n$, since if
$f\in X$, then $f(t)-f(t_1)=:g(t)\in X(T,H;t_1;0)$
(in the third case we use that $T\in \cal R$ and hence there exists
the embedding $H\hookrightarrow X$). For $n=1$ and $t_0\in T$
and $z_1=0$ we denote $X_0:=X_0(T,H):=X(T,H;t_0;0)$.
\par {\bf 6.2. Definitions.} We define a (non-Archimedean)
stochastic process $w(t,\omega )$ with values in $H$ 
as a random variable such that:
\par $(i)$ the differences $w(t_4,\omega )-w(t_3,\omega )$ 
and $w(t_2,\omega )-w(t_1,\omega )$ are independent
for each chosen $(t_1,t_2)$ and $(t_3,t_4)$ with $t_1\ne t_2$,
$t_3\ne t_4$, such that either $t_1$ or $t_2$ is not in the two-element set
$ \{ t_3,t_4 \} ,$ where $\omega \in \Omega ;$
\par $(ii)$ the random variable $\omega (t,\omega )-\omega (u,\omega )$ has 
a distribution $\mu ^{F_{t,u}},$ where $\mu $ is a probability 
$\bf K_s$-valued measure on $(X(T,H),{\sf B})$ from \S 6.1, 
$\mu ^g(A):=\mu (g^{-1}(A))$ for $g: X\to H$ such that
$g^{-1}({\cal R}_H)\subset \sf B$ and each $A\in
{\cal R}_H$, a continuous linear operator $F_{t,u}: X\to H$ is given by
the formula $F_{t,u}(w):=w(t,\omega )-w(u,\omega )$ 
for each $w\in L^q(\Omega ,{\sf F},\lambda ;X),$
where $1\le q\le \infty ,$ ${\cal R}_H$ is a covering ring of $H$ such that
$F_{t,u}^{-1}({\cal R}_H)\subset \sf B$ for each $t\ne u$ in $T$;
\par $(iii)$ we also put $w(0,\omega )=0,$  
that is, we consider a $\bf K$-linear subspace
$L^q(\Omega ,{\sf F},\lambda ;X_0)$
of $L^q(\Omega ,{\sf F},\lambda ;X)$,
where $\Omega \ne \emptyset $, $X_0$
is the closed subspace of $X$ as in \S 6.1.
\par {\bf 7. Definition.} Let $B$ and $q$ be as in
\S 2.1 and denote by $\mu _{q,B,\gamma }$ the corresponding $q$-Gaussian
$\bf K_s$-valued measure on $H$.
Let $\xi $ be a stochastic process with a real time
$t\in T\subset \bf R$ (see Definition 6.2),
then it is called a non-Archimedean $q$-Wiener process
with real time (and controlled by $\bf K_s$-valued measure), if
\par $(ii)'$ the random variable $\xi (t,\omega )-\xi (u,\omega )$ has
a distribution $\mu _{q,(t-u)B,\gamma }$ for each $t\ne u\in T$.
\par Let $\xi $ be a stochastic process with a non-Archimedean time
$t\in T\subset \bf F$, where $\bf F$ is a local field,
then $\xi $ is called a non-Archimedean $q$-Wiener process
with $\bf F$-time (and controlled by $\bf K_s$-valued measure), if
\par $(ii)"$ the random variable $\xi (t,\omega )-\xi (u,\omega )$ has
a distribution $\mu _{q,ln [\chi _{\bf F}(t-u)]B,\gamma }$
for each $t\ne u\in T$,
where $\chi _{\bf F}: {\bf F}\to \bf T$ is a continuous character
of $\bf F$ as the additive group (see \S 2.5 \cite{lulapm,lujms112}).
\par {\bf 8. Proposition.} {\it For each given $q$-Gaussian measure
a non-Archimedean $q$-Wiener process with real ($\bf F$ respectively)
time exists.}
\par {\bf Proof.} In view of Proposition 6
for each $t>u>b$ a random variable $\xi (t,\omega )-\xi (b,\omega )$
has a distribution $\mu _{q,(t-b)B,\gamma }$ for real time parameter.
If $t$, $u$, $b$ are pairwise different points in $\bf F$,
then $\xi (t,\omega )-\xi (b,\omega )$ has a distribution
$\mu _{q,ln [\chi _{\bf F}(t-b)]B,\gamma }$, since $ln [\chi _{\bf F}(t-u)]+
ln [\chi _{\bf F}(u-b)]=ln [\chi _{\bf F}(t-b)]$.
This induces the Markov quasimeasure $\mu ^{(q)}_{x_0,\tau }$
on $(\prod_{t\in T}(H_t,{\sf U}_t)),$ where
$H_t=H$ and ${\sf U}_t=Bco (H)$ for each $t\in T$.
In view of Theorem 2.39 \cite{lulapm,lujms112}
there exists an abstract probability space
$(\Omega ,{\sf F},\lambda )$, consequently, the corresponding
space $L(\Omega ,{\sf F},\lambda ,{\bf K_s})$ exists.
\par {\bf 9. Proposition.} {\it Let $\xi $ be a $q$-Gaussian process
with values in a Banach space $H=c_0(\alpha ,{\bf K})$
a time parameter $t\in T$ (controlled by a $\bf K_s$-valued measure)
and a positive definite correlation operator $B$ of trace class
and $\gamma =0$, where $card (\alpha )\le \aleph _0$,
either $T\subset \bf R$ or $T\subset \bf F$. Then either
$$(i)\quad \lim_{N\in \alpha }M_t [v^s_q(e^1(\xi (t,\omega ))^2+ ... +
v^s_q(e^N(\xi (t,\omega )))^2]=t Tr(B) \mbox{ or }$$
$$(ii)\quad \lim_{N\in \alpha }M_t [v^s_q(e^1(\xi (t,\omega ))^2+ ... +
v^s_q(e^N(\xi (t,\omega ))^2]=[ln (\chi _{\bf F}(t))] Tr(B)
\mbox{ respectively}.$$ }
\par {\bf Proof.} Define $\bf U_s$-valued moments \\
$m^q_k(e^{j_1},...,e^{j_k}):=
\int_Hv^s_{2q}(e^{j_1}(x))...v^s_{2q}(e^{j_k}(x))\mu _{q,B,\gamma }(dx)$ \\
for linear continuous functionals $e^{j_1},...,e^{j_k}$
on $H$ such that $e^l(e_j)=\delta ^l_j$, where
$\{ e_j: j\in \alpha \} $ is the standard orthonormal base in $H$.
\par Consider the operator
\par $(iii)$ $\mbox{ }_P\partial ^u\psi (x):=
F^{-1}({\hat f}_{u-1}(y) {\hat {\psi }}(y) |y|_p)(x)$, \\
where $f_u(x):=s^{-(1+u)\times ord_p(x)}/ \Gamma ^{{\bf K},s}(1+u)$ and
$F(f_u)(y) =\Gamma ^{{\bf K},s}(1+u)f_u(-y)|y|_p^{-1}$ (see \S 4),
where $F$ denotes the Fourier-Stieltjes operator defined
with the help of the $\bf K_s$-valued Haar measure $w$
on $Bco ({\bf K})$, $F (\psi )=:{\hat {\psi }}$, $Re (u)\ne -1$,
$\psi : {\bf K}\to \bf K_s$. Then
\par $(iv)$ $\mbox{ }_P\partial ^uf_b(x)=F^{-1}(
\Gamma ^{{\bf K},s}(u)f_{u-1}(-y)
\Gamma ^{{\bf K},s}(1+b)f_b(-y)|y|_p^{-1})$  $=f_{(u+b)}(x)$ \\
for each $u$ with $Re (u)\ne 0$, since  \\
$F^{-1}(s^{-(1+u+b)\times ord_p(-y)}|y|_p^{-1})(x)=
(\Gamma ^{{\bf K},s}(1+u+b))^{-1}s^{-(1+u+b)\times ord_p(-y)}(x)$. \\
For $u=1$ we write shortly
$\mbox{ }_P\partial ^1=\mbox{ }_P\partial $ and
$\mbox{ }_P\partial ^u_j$ means the operator of partial pseudo-differential
(with weight multiplier) given by Equation $(iii)$
by the variable $x_j$.
A function $\psi $ for which $\mbox{ }_P\partial ^u_j\psi $ exists
is called pseudo-differentiable (with weight multiplier)
of order $u$ by variable $x_j$. Then \\
$m^{q/2}_{2k}(e^{j_1},...,e^{j_{2k}})(\Gamma ^{{\bf K},s}(q/2))^{2k}$ 
$:= \int_Hs^{-q~ ord_p(x_{j_1})/2}...s^{-q~ ord_p(x_{j_{2k}})/2}
\mu _{q,B,\gamma }(dx)$ \\
$=\mbox{ }_P\partial ^{q/2}_{j_1}...\mbox{ }_P\partial ^{q/2}_{j_{2k}}
{\hat {\mu }}_{q,B,\gamma }(0)$ 
$=([\mbox{ }_PD^{q/2}]^{2k}{\hat {\mu }}(x))|_{x=0}.
(e^{j_1},...,e^{j_{2k}})$,  \\
where $(\mbox{ }_PD^{q/2}f(x)).e^j:=\mbox{ }_P\partial _jf(x)$. Therefore,
\par $(v)$ $m^{q/2}_{2k}(e^{j_1},...,e^{j_{2k}})
(\Gamma ^{{\bf K},s}(q/2))^{2k}$ \\ 
$=(k!)^{-1} [\mbox{ }_PD^{q/2}]^{2k}[B(v^s_q(z),v^s_q(z)]^k.
(e_{j_1},...,e_{j_{2k}})$
\par $=(k!)^{-1} \sum_{\sigma \in \Sigma _{2k}}
B_{\sigma (j_1),\sigma (j_2)}...B_{\sigma (j_{2k-1}),\sigma (j_{2k})}$, \\
since $\gamma =0$ and $\chi _{\gamma }(z)=1$,
where $\Sigma _k$ is the symmetric group of
all bijective mappings $\sigma $ of the set $\{ 1,...,k \} $
onto itself, $B_{l,j}:=B(e_j,e_l)$, since $Y^*=Y$ for
$Y=l_2(\alpha ,{\bf R})$.
Therefore, for each $B\in L_1$ and $A\in L_{\infty }$
we have $\int_HA(v_q(x),v_q(x))\mu _{q,B,0}(dx)=$
$\lim_{N\in \alpha }\sum_{j=1}^N\sum_{k=1}^N
A_{j,k}m^{q/2}_2(e_j,e_k)$ $=Tr(AB)$, since
${\bf C_s}\subset \bf U_s$ and algebraically $\bf C_s$ is
isomorphic with $\bf C$.
\par In particular for $A=I$ and $\mu _{q,tB,0}$ corresponding
to the transition measure of $\xi (t,\omega )$
we get Formula $(i)$ for a real time parameter,
using $\mu _{q,ln [\chi _{\bf F}(t)]B,0}$
we get Formula $(ii)$ for a time parameter belonging to $\bf F$,
since $\xi (t_0,\omega )=0$ for each $\omega $.
\par {\bf 10. Corollary.} {\it Let $H=\bf K$
and $\xi $, $B=1$, $\gamma $ be as in Proposition 9, then
$$(i)\quad M(\int_{t\in [a,b]} \phi (t,\omega )v^s_{2q}(d\xi (t,\omega ))=
M[\int_a^b\phi (t,\omega )dt]$$
for each $a<b\in T$ with real time, where
$\phi (t,\omega )\in L(\Omega ,{\sf U},\lambda ,C^0_0(T,{\bf R}))$
$\xi \in L(\Omega ,{\sf U},\lambda ,X_0(T,{\bf K}))$,
$(\Omega ,{\sf U},\lambda )$ is a probability measure space.}
\par {\bf Proof.} Since $\int_{t\in [a,b]}\phi (t,\omega )
v^s_{2q}(d\xi (t,\omega ))$ \\
$=\lim_{\max_j (t_{j+1}-t_j)\to 0}
\sum_{j=1}^N\phi (t_j,\omega ) v^s_q(\xi (t_{j+1},\omega )-
\xi (t_j,\omega ))$ for $\lambda $-almost all $\omega
\in \Omega $, since ${\bf C_s}\subset \bf U_s$ and $\bf C_s$ is
algebraically isomorphic with $\bf C$, then from the application
of Formula $9.(i)$ to each
$v^s_{2q} (\xi (t_{j+1},\omega )-\xi (t_j,\omega ))$
and the existence of the limit by finite partitions
$a=t_1<t_2<...<t_{N+1}=b$ of the segment $[a,b]$ it follows
Formula $10.(i)$.
\par {\bf 11. Definitions and Notes.} Consider a pseudo-differential operator
on $H=c_0(\alpha ,{\bf K})$ such that
$$(i)\quad {\sf A}=\sum_{0\le k\in {\bf Z};
j_1,...,j_k\in \alpha }(-i)^kb^k_{j_1,...,j_k}
\mbox{ }_P\partial _{j_1}...\mbox{ }_P\partial _{j_k},$$
where $b^k_{j_1,...,j_k}\in \bf R$,
$\mbox{ }_P\partial _{j_k}:=\mbox{ }_P\partial ^1_{j_k}$.
If there exists $n:=\max \{ k: b^k_{j_1,...,j_k}
\ne 0, j_1,...,j_k\in \alpha  \} $,
then $n$ is called an order of ${\sf A}$, $Ord ({\sf A})$,
where $\mbox{ }_P\partial _j$ is defined by Formula $9.(iii)$.
If ${\sf A}=0$, then by definition $Ord ({\sf A})=0$.
If there is not any such finite $n$, then $Ord ({\sf A})=\infty $.
We suppose that the corresponding form $\tilde A$ on $\bigoplus_kY^k$
is continuous into $\bf C$, where
$$(ii)\quad {\tilde A}(y)=-\sum_{0\le k\in {\bf Z};
j_1,...,j_k\in \alpha }(-i)^kb^k_{j_1,...,j_k}y_{j_1}...y_{j_k}/ ln s,$$
$y\in l_2(\alpha ,{\bf R})=:Y$.
If ${\tilde A}(y)>0$ for each $y\ne 0$ in $Y$, then
$\sf A$ is called strictly elliptic pseudodifferential operator.
\par Let $X$ be a complete locally $\bf K$-convex space,
let $Z$ be a complete locally $\bf U_s$-convex space.
For $0\le n\in \bf R$ a space of all functions $f: X\to Z$
such that $f(x)$ and $(\mbox{ }_PD^kf(x)).(y^1,...,y^{l(k)})$ are continuous
functions on $X$ for each $y^1,...,y^{l(k)}\in \{ e^1,e^2,e^3,... \}
\subset X^*$, $l(k):=[k]+sign \{ k \} $ for each $k\in
\bf N$ such that $k\le [n]$ and also for $k=n$
is denoted by $\mbox{ }_P{\cal C}^n(X,Z)$ and $f\in \mbox{ }_P
{\cal C}^n(X,Z)$ is called $n$ times continuously pseudodifferentiable, where
$[n]\le n$ is an integer part of $n$, $1> \{ n \} :=n-[n]\ge 0$
is a fractional part of $n$.
Then $\mbox{ }_P{\cal C}^{\infty }(X,Z):=\bigcap_{n=1}^{\infty }
\mbox{ }_P{\cal C}^n(X,Z)$
denotes a space of all infinitely pseudo-differentiable functions.
\par Embed $\bf R$ into $\bf C_s$ and consider the function
$v_2^s: {\bf U_p}\to {\bf R}\subset \bf C_s$, then for $t=v_2^s(\theta )$,
$\theta \in {\bf K}\subset \bf U_p$, put
$\partial _tu(t,x):=\lim_{\theta , {\bf K}, \theta \in {\bf K},
v_2^s(\theta )\to t} \mbox{ }_P\partial _{\theta }
u(v_2^s(\theta ),x)$ for $t\ge 0$, when it exists by the filter of
local subfields $\bf K$ in $\bf C_p$, which is correct, since
$v_2^s({\bf U_p})=[0,\infty )$, $\bigcup_{{\bf K}\subset \bf C_p}
{\bf K}$ is dense in $\bf C_p$, $\Gamma _{\bf C_p}=(0,\infty )\cap \bf Q$.
\par {\bf 12. Theorem.} {\it Let $\sf A$ be a
strictly elliptic pseudodifferential
operator on $H=c_0(\alpha ,{\bf K})$, $card (\alpha )\le \aleph _0$,
and let $t\in T=[0,b]\subset \bf R$. Suppose also that
$u_0(x-y)\in L(H,Bco (H),\mu _{t{\tilde A}},{\bf U_s})$
for each marked $y\in H$ as a function by $x\in H$,
$u_0(x)\in \mbox{ }_P{\cal C}^{Ord ({\sf A})}(H,{\bf U_s})$.
Then the non-Archimedean analog of the Cauchy problem
$$(i)\quad \partial _tu(t,x)={\sf A}u,\quad u(0,x)=u_0(x)$$
has a solution given by
$$(ii)\quad u(t,x)=\int_Hu_0(x-y)\mu _{t{\tilde A}}(dy),$$
where $\mu _{t{\tilde A}}$ is a $\bf K_s$-valued measure on $H$
with a characteristic functional ${\hat \mu }_{t{\tilde A}}(z)
:=s^{t{\tilde A}(v^s_2(z))}$.}
\par {\bf Proof.} In accordance with \S \S 2 and 11
we have $Y=l_2(\alpha ,{\bf R})$. The
function $s^{t{\tilde A}(v^s_2(z))}$
is continuous on $H\hookrightarrow H^*$ for each $t\in \bf R$
such that the family $H$ of continuous $\bf K$-linear functionals
on $H$ separates points in $H$. In view of Theorem
2.30 above it defines a tight measure on $H$ for each $t>0$.
The functional $\tilde A$ on each ball of radius $0<R<\infty $
in $Y$ is a uniform limit of its restrictions
${\tilde A}|_{\bigoplus_k[span_{\bf K}(e_1,...,e_n)]^k},$
when $n$ tends to the infinity, since $\tilde A$ is
continuous on $\bigoplus_kY^k$.
Since $u_0(x-y)\in L(H,Bco (H),\mu _{t{\tilde A}},{\bf U_s})$
and a space of cylindrical functions is dense in the latter Banach space
over $\bf U_s$,
then in view of Theorems $9.14, 9.21$ \cite{roo} and the Fubini theorem
it follows that $\lim_{P\to I}{\sf F}_{Px}u_0(Px)){\hat \mu }_{t\tilde A}
(y+Px)$ converges in $L(H,Bco (H),\mu _{t{\tilde A}},{\bf U_s})$
for each $t$, since $\mu _{t_1\tilde A}*\mu _{t_2\tilde A}=\mu_{(t_1+t_2)
\tilde A}$ for each $t_1$, $t_2$ and $t_1+t_2\in T$, where
$P$ is a projection on a finite dimensional over $\bf K$ subspace
$H_P:=P(H)$ in $H$, $H_P\hookrightarrow H$,
$P$ tends to the unit operator $I$ in the strong operator
topology, $F_{Px}u_0(Px)$ denotes a Fourier transform by
the variable $Px\in H_P$.
Consider a function $v:=F_x(u)$, then $\partial _tv(t,x)
=-{\tilde A}(v^s_2(x))v(t,x) ln s$, consequently,
$v(t,x)=v_0(x)s^{t{\tilde A}(v^s_2(x))}$. From $u(t,x)=F_x^{-1}
(v(t,x))$, where 
$F_x(u(t,x))=\lim_{n\to \infty }F_{x_1,...,x_n}u(t,x)$.
Therefore, $u(t,x)=u_0(x)*[F_x^{-1}({\hat \mu }_{t{\tilde A}})]=$
$\int_Hu_0(x-y)\mu _{t{\tilde A}}(dy)$,
since $u_0(x-y)\in L(H,Bco (H),\mu _{t{\tilde A}},{\bf U_s})$
and $\mu _{t{\tilde A}}$ is the tight measure on $Bco (H)$.
\par {\bf 13. Note.} In the particular case of $Ord ({\sf A})=2$
and ${\tilde A}$ corresponding to the Laplace operator, that is,
${\tilde A}(y)=\sum_{l,j} g_{l,j}y_ly_j$, Equation
$12.(i)$ is (the non-Archimedean analog of) the heat equation
on $H$.
\par For $Ord ({\sf A})<\infty $ the form ${\tilde A}_0(y)$ corresponding
to sum of terms with $k=Ord ({\sf A})$ in Formula $11.(ii)$ is called
the principal symbol of operator $\sf A$. If ${\tilde A}_0(y)>0$
for each $y\ne 0$, then $\sf A$ is called an elliptic pseudodifferential
operator. Evidently, Theorem 12 is true for elliptic $\sf A$
of $Ord ({\sf A})<\infty $.
\par {\bf 14. Remark and Definitions.} Let linear spaces $X$ over $\bf K$
and $Y$ over $\bf R$ be as in \S 4 and $B$ be a symmetric
nonnegative definite (bilinear) operator on a dense $\bf R$-linear subspace
${\sf D}_{B,Y}$ in $Y^*$. A quasi-measure $\mu $ with a characteristic
functional
$${\hat \mu }(\zeta ,x):=s^{\zeta B(v^s_q(z),v^s_q(z))}\chi _{\gamma }(z)$$
for a parameter $\zeta \in \bf C_s$ with $Re (\zeta )\ge 0$
defined on ${\sf D}_{q,B,X}$
we call an $\bf U_s$-valued (non-Archimedean analog of Feynman)
quasi-measure and we denote it
by $\mu _{q,\zeta B,\gamma }$ also, where ${\sf D}_{q,B,X}:=
\{ z\in X^*:$ $\mbox{there exists}$ $j\in \Upsilon $
$\mbox{such that}$ $z(x)=z_j(\phi _j(x))$ $\forall x\in X,$
$v^s_q(z)\in {\sf D}_{B,Y} \} $.
\par {\bf 15. Proposition.} {\it Let $X={\sf D}_{q,B,X}$
and $B$ be positive definite, then
for each function $f(z):=\int_X\chi _z(x)\nu (dx)$ with
an $\bf U_s$-valued tight measure $\nu $ of finite norm
and each $Re (\zeta )>0$ there exists
$$(i)\quad \int_Xf(z)\mu _{\zeta B}(dz)=
\lim_{P\to I}\int_Xf(Pz)\mu ^{(P)}_{\zeta B}(dz)$$
$$=\int_Xs^{(\zeta B(v_q(z),v_q(z)))}\chi _{\gamma }(z)\nu (dz),$$
where $\mu ^{(P)}(P^{-1}(A)):=\mu (P^{-1}(A))$ for each
$A\in Bco (X_P)$, $P: X\to X_P$ is a projection on a
$\bf K$-linear subspace $X_P$, a convergence $P\to I$ is
considered relative to a strong operator topology.}
\par {\bf Proof.} From the use of the projective limit decomposition
of $X$ and Theorem 2.37 \cite{lulapm,lujms112}
it follows, that there exists \\
$(ii)\quad \int_Xf(z)\mu _{\zeta B}(dz)=
\lim_{P\to I}\int_Xf(Pz)\mu ^{(P)}_{\zeta B}(dz)$.
Then for each finite dimensional over $\bf K$ subspace $X_P$ \\
$(iii)\quad \int_Xf(Pz)\mu ^{(P)}_{\zeta B}(dz)=\int_{X_P}
\{ s^{\zeta B(v^s_q(z),v^s_q(z)))} \chi _{\gamma }(z) \}
|_{X_P}\nu ^{X_P}(dz),$   \\
since $\nu $ is tight and hence each $\nu ^{X_P}$ is tight.
Each measure $\nu _j$ is tight,
then due to Lemma 2.3 and \S 2.5 \cite{lulapm,lujms112}
there exists the limit \\
$\lim_{P\to I} \int_{X_P}\{ s^{\zeta B(v^s_q(z),v^s_q(z))}
\chi _{\gamma }(z) \} |_{X_P} \nu ^{X_P}(dz)$\\
$= \int_Xs^{\zeta B(v^s_q(z),v^s_q(z)))}\chi _{\gamma }(z)\nu (dz).$
\par {\bf 16. Proposition.} {\it If conditions of Proposition
15 are satisfied and
$$(i)\quad f(Px)\in L(X_P,Bco (w^{X_P}),{\bf U_s}) $$
for each finite dimensional over $\bf K$ subspace $X_P$ in $X$ and
$$(ii)\quad \lim_{R\to \infty }\sup_{|x|\le R} |f(x)|=0,$$
then Formula $15.(i)$ is accomplished for $\zeta $ with
$Re (\zeta )=0$, where $w^{X_P}$ is a nondegenerate $\bf K_s$-valued
Haar measure on $X_P$.}
\par {\bf Proof.} In view of Theorem 2.37 \cite{lulapm,lujms112}
for the consistent
family of measures $\{ f(Px)\mu ^{X_P}_{q,iB,\gamma }(d Px): P \} $
(see \S 2.36 \cite{lulapm,lujms112})
there exists a measure on $(X,{\cal R})$, where
projection operators $P$ are associated with a chosen basis in $X$.
The finite dimensional over $\bf K$ distribution
$\mu ^{X_P}_{q,iB,\gamma }/w^{X_P}(dx)=F^{-1}(
{\hat \mu }_{q,iB,\gamma })|_{X_P})$ is in $C_{\infty }(X_P,{\bf U_s})$
due to Theorem $9.21$ \cite{roo}, since ${\hat {\mu }}\in L(X_P,Bco (X_p),
w^{X_P},{\bf U_s})$. In view of Condition $16.(i,ii)$ above
and the Fubini theorem and using the Fourier-Stieltjes
transform we get Formulas $15.(ii,iii)$. From the taking the limit
by $P\to I$ Formula $15.(i)$ follows. This means that
$\mu _{q,\zeta B,\gamma }$ exists in the sence of distributions.
\par {\bf 17. Remark.} Put
$$(i)\quad \mbox{ }_F\int_X f(x)\mu _{q,iB,\gamma }(dx)
:=\lim_{\zeta \to i} \int_Xf(x)\mu _{q,\zeta B,\gamma }(dx)$$
if such limit exists. If conditions of Proposition 16 are
satisfied, then $\psi (\zeta ):=\int_Xf(x)\mu _{q,\zeta B,\gamma }(dx)$
is the pseudo-differentiable of order $1$ function by $\zeta $
on the set $\{ \zeta \in {\bf C_s}: Re (\zeta )>0 \} $
and it is continuous on the subset $\{ \zeta \in {\bf C_s}:
Re (\zeta )\ge 0 \} $, consequently,
$$(ii)\quad \mbox{ }_F\int_X f(x)\mu_{q,iB,\gamma }(dx)=
\int_Xs^{ \{ i B(v^s_q(x),v^s_q(x)) \} } \chi _{\gamma }(x)\nu (dx).$$
\par Above non-Archimedean analogs of Gaussian measures with specific
properties were defined. Nevertheless, there do not exist usual
Gaussian $\bf K_s$-valued measures on non-Archimedean Banach spaces.
\par {\bf 18. Theorem.} {\it Let $X$ be a Banach space of separable type
over a locally compact non-Archimedean field $\bf K$.
Then on $Bco (X)$ there does not exist a nontrivial
$\bf K_s$-valued (probability) usual Gaussian measure.}
\par {\bf Proof.} Let $\mu $ be a nontrivial usual Gaussian $\bf K_s$-valued
measure on $Bco (X)$. Then by the definition its characteristic
functional ${\hat {\mu }}$ must be satisfying Conditions $2.5.(3,5)$
\cite{lulapm,lujms112}
$\bf U_s$-valued function and $\lim_{|y|\to \infty }
{\hat {\mu }}(y)=0$ for each $y\in X^* \setminus \{ 0 \} $,
where $X^*$ is the topological conjugate space to $X$
of all continuous $\bf K$-linear functionals $f: X\to \bf K$.
Moreover, there exist a $\bf K$-bilinear functional $g$ and
a compact nondgenerate $\bf K$-linear operator $T: X^*\to X^*$
with $ker (T)= \{ 0 \} $ and a marked vector $x_0\in X$
such that ${\hat {\mu }}_{x_0}(y)=f(g(Ty,Ty))$ for each
$y\in X^*$, where $\mu _{x_0}(dx):=\mu (-x_0+dx)$, $x\in X$.
Since $\bf K$ is locally compact, then $X^*$ is nontrivial
and separates points of $X$ (see \cite{nari,roo}).
Each one-dimensional over $\bf K$ projection of a Gaussian
measure is a Gaussian measure and products of Gaussian measures
are Gaussian measures, hence convolutions of Gaussian measures
are also Gaussian measures. Therefore, ${\hat {\mu }}_{x_0}:
X^*\to \bf U_s$ is a nontrivial character: ${\hat {\mu }}_{x_0}(y_1+y_2)=
{\hat {\mu }}_{x_0}(y_1) {\hat {\mu }}_{x_0}(y_2)$
for each $y_1$ and $y_2$ in $X^*$. If $char ({\bf K})=0$
and $\bf K$ is a non-Archimedean field, then
there exists a prime number $p$ such that $\bf Q_p$
is the subfield of $\bf K$. Then ${\hat {\mu }}(p^ny)=
({\hat {\mu }}(y))^{p^n}$ for each $n\in \bf Z$ and $y\in
X^*\setminus \{ 0 \} $, particularly, for $n\in \bf N$ tending
to the infinity we have $\lim_{n\to \infty } p^ny = 0$
and $\lim_{n\to \infty } {\hat {\mu }}_{x_0}(p^ny)=1$,
$\lim_{n\to \infty } {\hat {\mu }}_{x_0}(y))^{p^n}=0$, since
$s\ne p$ are primes, $\lim_{n\to \infty }{\hat {\mu }}_{x_0}(p^{-n}y)=0$
and $|{\hat {\mu }}_{x_0}(y)|<1$ for $y\ne 0$.
This gives the contardiction, hence $\bf K$ can not be
a non-Archimedean field of zero characteristic.
Suppose that $\bf K$ is a non-Archimedean field of characteristic
$char ({\bf K})=p>0$, then $\bf K$ is isomorphic with the field
of formal power series in variable $t$ over a finite field $\bf F_p$.
Therefore, ${\hat {\mu }}_{x_0}(py)=1$, but ${\hat {\mu }}_{x_0}(y)^p\ne 1$
for $y\ne 0$, since $\lim_{n\to \infty } {\hat {\mu }}_{x_0}
(t^{-n}y)=0$. This contradicts the fact that ${\hat {\mu }}_{x_0}$
need to be the nontrivial character, consequently, $\bf K$ can
not be a non-Archimedean field of nonzero characteristic as well.
It remains the classical case of $X$ over $\bf R$ or $\bf C$,
but the latter case reduces to $X$ over $\bf R$ with the help
of the isomorphism of $\bf C$ as the $\bf R$-linear space with $\bf R^2$.
\par {\bf 19. Theorem.} {\it Let $\mu _{q,B,\gamma }$ and
$\mu _{q,B,\delta }$ be two $q$-Gaussian $\bf K_s$-valued measures.
Then $\mu _{q,B,\gamma }$ is equivalent to $\mu _{q,B,\delta }$
or $\mu _{q,B,\gamma }\perp \mu _{q,B,\delta }$ according to
$v^s_q(\gamma -\delta ) \in B^{1/2}({\sf D}_{B,Y})$ or not.
The measure $\mu _{q,B,\gamma }$ is orthogonal to $\mu _{g,B,\delta }$,
when $q\ne g$. Two measures $\mu _{q,B,\gamma }$
and $\mu _{g,A,\delta }$ with positive definite nondegenerate
$A$ and $B$ are either equivalent or orthogonal.}
\par {\bf 20. Theorem.} {\it The measures $\mu _{q,B,\gamma }$
and $\mu _{q,A,\gamma }$ are equivalent if and only if there exists
a positive definite bounded invertible operator $T$ such that
$A=B^{1/2}TB^{1/2}$ and $T-I\in L_2(Y^*)$. }
\par {\bf Proof.} Using the projective limit reduce consideration to the
Banach space $X$. Let $z\in X$ be a marked vector and $P_z$ be a
projection operator on $z\bf K$ such that $P_z^2=P_z$,
$z=\sum_jz_je_j$, then the characteristic functional of the projection
$\mu _{q,B,\gamma }^{z\bf K}$ of $\mu _{q,B,\gamma }$ has the form
${\hat {\mu }}_{q,B,\gamma }^{z\bf K}=s^{[(\sum_{i,j}B_{i,j}v_q^s(z_i)
v_q^s(z_j))v_{2q}^s(\xi )]}\chi _{\gamma (z)}(\xi )$
for each vector $x=\xi z$, where
each $z_j$ and $\xi \in \bf K$, since $v_{2q}^s(\xi )=(v_q^s(\xi ))^2$.
Choose a sequence $ \{ \mbox{ }_nz: n \} $ in $X$ such that
it is the orthonormal basis in $X$ and the operator $G: X\to X$
such that $G\mbox{ }_nz=\mbox{ }_na \mbox{ }_nz$ with $\mbox{ }_na\ne 0$
for each $n\in \bf N$ and there exists $G^{-1}: G(X)\to X$
such that it induces the operator $C$ on a dense subspace ${\cal D}(Y)$
in $Y$ such that $CBC: Y\to Y$ is invertible and
$\| CBC \| $ and $\| (CBC)^{-1} \| \in [|\pi |, |\pi |^{-1}]$. Then \\
$\mu _{q,A,\gamma }(dx)/ \mu _{q,B,\gamma }(dx)=
\lim_{n\to \infty } [\mu _{q,A,\gamma }^{V_n}(dx^n)/
\lambda ^{V_n}(dx^n)] [\mu _{q,B,\gamma }^{V_n}(dx^n)/
\lambda ^{V_n}(dx^n)]^{-1}$, where $V_n:=span_{\bf K}(\mbox{ }_jz:
j=1,...,n)$, $x_n\in V_n$. Consider $x_n=G^{-1} (y_n)$,
where $y_n\in G (V_n)$, then \\
$[\mu _{q,B,\gamma }^{V_n}(G^{-1}dy^n)/
\lambda ^{V_n}(G^{-1}dy^n)]$ and
$[\mu _{q,B,\gamma }^{V_n}(G^{-1}dy^n)/ \lambda ^{V_n}(G^{-1}dy^n)]^{-1}$ \\
are in $L(\lambda ^{V_n}(G^{-1}dy^n))$ for each $n$
such that there exists $m\in \bf N$ for which \\
$\| [\mu _{q,B,\gamma }^{V_n} (G^{-1}dy^n)/
\lambda ^{V_n} (G^{-1}dy^n)] \| $ and
$\| [\mu _{q,B,\gamma }^{V_n} (G^{-1}dy^n)/
\lambda ^{V_n} (G^{-1}dy^n)]^{-1} \| \in [ |\pi |, |\pi |^{-1} ]$ \\
for each $n>m$, where $\| * \| $ is taken in
$L(\lambda ^{V_n}(G^{-1}dy^n ))$.
Then $N_{\mu _{q,CBC,\gamma G^{-1}}^{V_n}}(y^n)\in
[|\pi |, |\pi |^{-1}]$ for each $n>m$. Then the existence of
$\mu _{q,A,\gamma }(dx)/ \mu _{q,B,\gamma }(dx)\in
L(\mu _{q,B,\gamma })$ is provided by using operator $G$ and
the consideration of characteristic functionals of measures, Theorem $3.5$
and the fact that the Fourier-Stieltjes transform $F$ is
the isomorphism of Banach algebras $L({\bf K}, Bco ({\bf K}), v, {\bf U_s})$
with $C_{\infty }({\bf K},{\bf U_s})$, where $v$ denotes the Haar
normalized by $v(B({\bf K},0,1))=1$ $\bf K_s$-valued measure on $\bf K$.
If $g\ne q$ then the measure $\mu _{q,B,\gamma }$
is orthogonal to $\mu _{g,B,\delta }$, since \\
$\lim_{R>0, R+n\to \infty } \sup_{x\in X^c_{R,n}}
|(\mu _{q,B,\gamma })_{X_n}/(\mu _{g,B,\delta })_{X_n}|(x)=0$ \\
for each $q>g$ due to Formula $4.(ii)$,
where $X_n:=span_{\bf K}(e_m: m=n,n+1,...,2n)$,
$X^c_{R,n}:=X_n\setminus B(X_n,0,R)$ ,
$(\mu _{q,B,\gamma })_{X_n}$ is the projection of the measure
$\mu _{q,B,\gamma }$ on $X_n$. Each term $\beta _j$ in Theorem $3.5$ is
in $[0,1]\subset \bf R$, consequently, the product in this theorem
is either converging to a positive constant or diverges to zero,
hence two measures $\mu _{q,B,\gamma }$
and $\mu _{g,A,\delta }$ are either equivalent or orthogonal.
\par {\bf 21. Theorem.} {\it Let $X$ be a Banach space of separable
type over a locally compact non-Archimedean field $\bf K$
and $J$ be a dense proper $\bf K$-linear subspace in $X$ such
that the embedding operator $T: J\hookrightarrow X$ is compact
and nondegenerate, $ker (T)= \{ 0 \} $.
Then a set ${\cal M}(X,J)$ of probability $\bf K_s$-valued
measures $\mu $ on $Bco (X)$ quasi-invariant relative to $J$
is of cardinality ${card ({\bf K_s})}^{\sf c}$.
If $J'$, $J'\subset J$, is also a dense $\bf K$-linear subspace in $X$,
then ${\cal M}(X,J')\supset {\cal M}(X,J)$.}
\par {\bf Proof.} Since $X$ is of separable type over $\bf K$,
then we can choose for a given compact operator $T$ an orthonormal base
in $X$ in which $T$ is diagional and $X$ is isomorphic with
$c_0$ over $\bf K$ such that in its standard base
$ \{ e_j: j\in {\bf N} \} $ the operator $T$ has the form
$Te_j=a_je_j$, $0\ne a_j\in \bf K$ for each $j\in \bf N$,
$\lim_{j\to \infty }a_j=0$. As in Theorem 3.15 \cite{lulapm,lujms112} 
take $g_n\in L({\bf K},Bco ({\bf K}),w'(dx/a_n),{\bf K_s})$,
$g_n(x)\ne 0$ for $v$-a.e. $x\in \bf K$ and $\| g_n \| =1$ for each $n$,
for which converges $\prod_{n=1}^{\infty }
\beta _n>0$ for each $y\in J$ and such that
$\prod_{n=1}^mg_n(x_n)w'(dx_n/a_n)=:\nu _{L_n}(dx^n)$ satisfies
conditions of Lemma $2.3$ \cite{lulapm,lujms112},
where $\beta _n:=\| \rho _n \|_{\phi _n}$,
$0\ne a_n\in \bf K$ for each $n\in \bf N$,
$\rho _n(x) := \mu _n(dx)/\nu _n(dx)$, $\phi _n(x):=N_{\lambda _n}(x)$,
$\lambda _n(dx):= g_n(x) w'(dx/a_n)$, then use 
Theorem $3.5$ \cite{lulapm,lujms112}
for the measure $\nu _n(dx):=g_n(x)w'(dx/a_n)$ and
$\mu _n(dx):=\nu _n(-y_n+dx)$, $x^n:=(x_1,...,x_n)$,
$x_1,...,x_n\in \bf K$ for each $n\in \bf N$.
The family of such sequences of functions $\{ g_n: n\in {\bf N} \} $
has the cardinality ${card ({\bf K_s})}^{\sf c}$, since in
$L(\nu )$ the subspace of step functions is dense and $card (Bco (X))=\sf c$.
The family of all $\{ g_n: n \} $ satisfying conditions above
for $J$ also satisfies such conditions for $J'$.
From which the latter statement of this theorem follows.
\par {\bf 22. Theorem.} {\it Let $X$ be a Banach space of separable
type over a locally compact non-Archimedean field $\bf K$
and $J$ be a dense proper $\bf K$-linear subspace in $X$ such
that the embedding operator $T: J\hookrightarrow X$ is compact
and nondegenerate, $ker (T)= \{ 0 \} $, $b\in \bf C$.
Then a set ${\cal P}_b(X,J)$ of probability
$\bf K_s$-valued measures $\mu $ on $Bco (X)$
quasi-invariant and pseudo-differentiable of order $b$
relative to $J$ is of cardinality $card ({\bf K_s})^{\sf c}$.
If $J'$, $J'\subset J$, is also a dense $\bf K$-linear subspace in $X$,
then ${\cal P}_b(X,J')\supset {\cal P}_b(X,J)$.}
\par {\bf Proof.} As in \S 21 choose for $T$ an orthonormal base
in $X$ in which $T$ is diagional and $X$ is isomorphic with
$c_0$ over $\bf K$ such that in its standard base
$\{ e_j: j\in {\bf N} \} $  the operator $T$ is characterized by
$Te_j=a_je_j$, $0\ne a_j\in \bf K$ for each $j\in \bf N$,
$\lim_{j\to \infty }a_j=0$. Take $g_n$ from \S 21,
where $g_n\in L({\bf K},Bf ({\bf K}),w'(dx/a_n),{\bf K_s})$,
satisfy conditions there and such that there exists
$\lim_{m\to \infty } PD(b,\prod_{n=1}^mg_n(xz))\in L(X, Bco (X),
\nu ,{\bf F})$ by the variable $x$ for each $z\in J$, where $x\in \bf K$,
${\bf K_s}\cup {\bf C_s}\subset \bf F$, $\bf F$ is a non-Archimedean field.
Evidently, ${\cal P}_b(X,J)\subset {\cal M}(X,J)$.
The family of such sequences of functions $\{ g_n: n\in {\bf N} \} $
has the cardinality $card ({\bf K_s})^{\sf c}$, since in
$L (\nu )$ the subspace
of step functions is dense and the condition of pseudo-differentiability
is the integral convergence condition (see \S \S 4.1 and 4.2
\cite{lulapm,lujms112}).


\begin{thebibliography}{299}
\bibitem[AK91]{alkar} Albeverio, S., Karwowski, W.:
Diffusion on $p$-adic numbers, 86--99. In: Ito, K.,
Hida, T. (eds.). Gaussian random
fields. Nagoya 1990. World Scientific, River Edge, NJ (1991)
\bibitem[ADV88]{ardrvo}  Aref'eva, I.Ya., Dragovich, B., Volovich, I.V.:
On the $p$-adic summability of the anharmonic oscillator.
Phys. Lett., \textbf{ B 200}, 512--514 (1988)
\bibitem[BV97]{byvo}  Bikulov, A. H., Volovich, I.V.: 
$p$-Adic Brounian motion.
Izvest. Russ. Acad. Sci. Ser. Math., \textbf{ 61: 3}, 75--90  (1997)
\bibitem[Bou63-69]{boui}  Bourbaki, N.: Int\'egration. Livre VI. 
Fasc. XIII, XXI, XXIX, XXXV. Ch. 1--9. Hermann, Paris
(1965, 1967, 1963, 1969).
\bibitem[Cas02]{cas}  Castro, C.: Fractal strings as an alternative
justification for El Naschie's cantorian spacetime and the fine
structure constants. Chaos, Solitons and Fractals,
\textbf{ 14}, 1341--1351 (2002)
\bibitem[DF91]{dal} Dalecky, Yu.L., Fomin, S.V.: Measures and differential
equations in infinite-dimensional spaces. Kluwer Acad. Publ.,
Dordrecht (1991)
\bibitem[Dia84]{diar}  Diarra, B.: Ultraproduits ultrametriques
de corps values. Ann. Sci. Univ. Clermont II, S\'er. Math., 
\textbf{ 22}, 1--37 (1984)
\bibitem[DD00]{djdr} Djordjevi\'c, G.S., Dragovich, B.:
$p$-Adic and adelic harmonic oscillator with a time-dependent
frequency. Theor. and Math. Phys., \textbf{ 124: 2}, 1059--1067 (2000)
\bibitem[Eng86]{eng}  Engelking, R.: General topology. Mir, Moscow (1986)
\bibitem[Esc95]{esc} Escassut, A.: Analytic elements in $p$-adic analysis.
World Scientific, Singapore (1995)
\bibitem[Eva88]{evans}  Evans, S.N.: Continuity properties
of Gaussian stochastic processes indexed by a local field.
Proceed. Lond. Math. Soc. Ser. 3, \textbf{ 56}, 380--416 (1988)
\bibitem[Eva89]{evans2}  Evans, S.N.: Local field Gaussian
measures, 121--160. In: Cinlar, E., et.al. (eds.)
Seminar on Stochastic Processes 1988. Birkh\"auser, Boston (1989) 
\bibitem[Eva91]{evans3}  Evans, S.N.: Equivalence and
perpendicularity of local field Gaussian
measures, 173--181. In:  Cinlar, E., et.al. (eds.)
Seminar on Stochastic Processes 1990. Birkh\"auser, Boston (1991) 
\bibitem[Eva93]{evans4}  Evans, S.N.: Local field Brownian motion.
J. Theoret. Probab. \textbf{ 6}, 817--850 (1993)
\bibitem[GV61]{gevi} Gelfand, I.M., Vilenkin, N.Ya.: Some applications
of harmonic analysis. Generalized functions.
\textbf{ 4} Fiz.-Mat. Lit., Moscow (1961)
\bibitem[Jan98]{yojan} Jang, Y.: Non-Archimedean quantum mechanics.
Tohoku Math. Publ. N \textbf{ 10} (1998)
\bibitem[Khr90]{khrum} Khrennikov, A.Yu.: Mathematical methods of
non-Archimedean physics. Russ. Math. Surv., \textbf{ 45: 4}, 79--110 (1990)
\bibitem[Khr91]{khrif} Khrennikov, A.Yu.:
Generalized functions and Gaussian
path integrals. Russ. Acad. Sci. Izv. Mat.,
\textbf{ 55}, 780--814 (1991)
\bibitem[Khr99]{khipr} Khrennikov, A.: Interpretations of probability.
VSP, Utrecht (1999)
\bibitem[Kob77]{kobl} Koblitz, N.: $p$-adic numbers, $p$-adic analysis
and zeta functions. Springer-Verlag, New York, 1977.
\bibitem[Lud96]{lu1} Ludkovsky, S.V.:
Measures on groups of diffeomorphisms of
non-Archimedean Banach manifolds. Russ. Math. Surv.,
\textbf{ 51: 2}, 338--340 (1996)
\bibitem[Lud96c]{lulapm}  Ludkovsky, S.V.: Quasi-invariant and
pseudo-differentiable measures
on a non-Archimedean Banach space. I, II. Los Alamos Preprints
{\bf math.GM/0106169} and {\bf math.GM/0106170} (http://xxx.lanl.gov/;
earlier version: ICTP {\bf IC/96/210},  October 1996, 50 pages
http://www.ictp.trieste.it/;
VINITI [Russ. Inst. of Sci. and Techn. Inform.],
Deposited Document \textbf{ 3353-B97}, 78 pages
(17 November 1997))
\bibitem[Lud98b]{luseamb} Ludkovsky, S.V.:
Irreducible unitary representations
of non-Archimedean groups of diffeomorphisms. Southeast Asian 
Bull. of Math., \textbf{ 22}, 419--436 (1998)
\bibitem[Lud98s]{luumn983}  Ludkovsky, S.V.: Quasi-invariant measures on
non-Archimedean semigroups of loops. Russ. Math. Surv.,
\textbf{ 53: 3}, 633--634 (1998)
\bibitem[Lud99a]{lubp99}  Ludkovsky, S.V.:
Properties of quasi-invariant measures on
topological groups and associated algebras.
Annales Math. B. Pascal, \textbf{ 6: 1}, 33--45 (1999) 
\bibitem[Lud99s]{luumn995} Ludkovsky, S.V.:
Non-Archimedean polyhedral expansions
of ultrauniform spaces. Russ. Math. Surv.,
\textbf{ 54: 5}, 163--164 (1999)
(detailed version: Los Alamos National Laboratory, USA.
Preprint {\bf math.AT/0005205}, 39 pages, May 2000)
\bibitem[Lud99t]{lutmf99}  Ludkovsky, S.V.:
Measures on groups of diffeomorphisms
of non-Archimedean manifolds, representations of groups and their 
applications. Theoret. and Math. Phys., \textbf{ 119: 3}, 698--711 (1999)
\bibitem[Lud00a]{lubp2}   Ludkovsky, S.V.:
Quasi-invariant measures on non-Archimedean
groups and semigroups of loops and paths, their representations. I, II.
Annales Math. B. Pascal, \textbf{ 7: 2}, 19--53, 55--80 (2000) 
\bibitem[Lud00f]{lufpm}  Ludkovsky, S.V.:
Non-Archimedean polyhedral decompositions of
ultrauniform spaces. Fundam. i Prikl. Math.
\textbf{ 6: 2}, 455--475 (2000)
\bibitem[Lud01f]{lufpmsp} Ludkovsky, S.V.: Stochastic processes
on groups of diffeomorphisms and loops of real, complex and
non-Archimedean manifolds. Fundam. i Prikl. Math.,
\textbf{ 7: 4}, 1091--1105 (2001)
\bibitem[Lud01s]{luumn01}  Ludkovsky, S.V.:
Representations of topological groups generated by Poisson measures.
Russ. Math. Surv., \textbf{ 56: 1}, 169--170 (2001)
\bibitem[Lud02a]{luanma}  Ludkovsky, S.V.:
Quasi-invariant and pseudo-differentiable
real-valued measures on a non-Archimedean Banach space.
Analysis Math., \textbf{ 28}, 287--316 (2002)
\bibitem[Lud02b]{luseamb2}  Ludkovsky, S.V.:
Poisson measures for topological groups and their representations.
Southeast Asian Bull. Math., \textbf{ 25: 4}, 653--680 (2002)
\bibitem[Lud0321]{luijmms1}  Ludkovsky, S.V.:
Stochastic processes on non-Archimedean Banach spaces.
Int. J. of Math. and Math. Sci., \textbf{ 2003: 21}, 1341--1363 (2003)
\bibitem[Lud0341]{luijmms2}  Ludkovsky, S.V.:
Stochastic antiderivational equations on non-Archimedean
Banach spaces. Int. J. of Math. and Math. Sci.,
\textbf{ 2003: 41}, 2587--2602 (2003)
\bibitem[Lud0348]{luijmms3} Ludkovsky, S.V.:
Stochastic processes on totally disconnected topological groups.
Int. J. of Math. and Math. Sci., \textbf{ 2003: 48}, 3067--3089 (2003)
\bibitem[Lud02j]{lujms112} Ludkovsky, S.V.:
Quasi-invariant and pseudo-differentiable
measures on non-Archimedean Banach spaces with values
in non-Archimedean fields. J. Math. Sci., is accepted to publication,
54 pages
\bibitem[Lud03s2]{luumn582} Ludkovsky, S.V.:  
Quasi-invariant and pseudo-differentiable measures on 
non-Archimedean Banach spaces. Russ. Math. Surv.,
\textbf{ 58: 2}, 167--168 (2003)
\bibitem[LD02]{luddiaop} Ludkovsky, S., Diarra, B.:
Spectral integration and spectral theory for non-Archimedean 
Banach spaces. Int. J. Math. and Math. Sci.,
\textbf{ 31: 7}, 421--442 (2002)
\bibitem[LK02]{lukhr}  Ludkovsky, S.V., Khrennikov, A.:
Stochastic processes on non-Archimedean spaces with values in
non-Archimedean fields. Markov Processes and Related Fields,
\textbf{ 8}, 1--34 (2002)
\bibitem[NB85]{nari} Narici, L., Beckenstein, E.:
Topological vector spaces. Marcel Dekker Inc., New York (1985)
\bibitem[Roo78]{roo} Rooij, A.C.M. van.:
Non-Archimedean functional analysis. Marcel Dekker Inc., New York (1978)
\bibitem[Sat94]{sat94} Sato, T.: Wiener measure on certain
Banach spaces over non-Archimedean local fields. Compositio Math.
\textbf{ 93}, 81--108 (1994)
\bibitem[Sch84]{sch1}  Schikhof, W.H.: Ultrametric calculus.
Camb. Univ. Press, Cambridge (1984)
\bibitem[VV89]{vla2} Vladimirov, V.S., Volovich, I.V.
Comm. Math. Phys., \textbf{123 }, 659--676 (1989)
\bibitem[VVZ94]{vla3}  Vladimirov, V.S., Volovich, I.V., Zelenov, E.I.:
$p$-Adic analysis and mathematical physics. Fiz.-Mat. Lit., Moscow (1994)
\bibitem[Wei73]{wei}  Weil, A.: Basic number theory. Springer, Berlin
(1973)
\end{thebibliography}
\end{document}